\documentstyle[12pt,graphicx,amsfonts]{article}
\newcommand{\vs}{\vspace{0.3cm}}

\def\thebibliography#1
{\begin{center}{\normalsize\bf References}\end{center}
\list{[\arabic{enumi}]}%
{\settowidth\labelwidth{[#1]}\leftmargin \labelwidth
\advance\leftmargin\labelsep
\usecounter{enumi}}
\sloppy\clubpenalty4000\widowpenalty4000}
\begin{document}

\makeatletter

%%% local diagrams
\def \minuscrossing{\raisebox{-4pt}{\begin{picture}(28,16)(-2,0)  
\thicklines
\put(25,-1){\vector(-3,2){26}} \put(15,10){\vector(3,2){10}}
\put(-1,-1){\line(3,2){10}} \end{picture}}}
\def \pluscrossing{\raisebox{-4pt}{\begin{picture}(28,16)(-2,0)  
\thicklines
\put(-1,-1){\vector(3,2){26}} \put(15,6){\line(3,-2){10}}
\put(10,10){\vector(-3,2){10}} \end{picture}}}
\def \opencrossing{\raisebox{-4pt}{\begin{picture}(28,16)(-2,0)  
\thicklines
\put(10,11){\vector(-2,1){10}} \put(10,11){\line(0,-1){6}}
\put(14,11){\vector(2,1){10}} \put(0,0){\line(2,1){10}}
\put(14,5){\line(0,1){6}}
\put(14,5){\line(2,-1){10}} \end{picture}}}
\def \doublepoint{\raisebox{-4pt}{\begin{picture}(28,16)(-2,0)  
\thicklines
\put(0,0){\vector(3,2){24}} \put(24,0){\vector(-3,2){24}}
\put(12,8){\circle*{3}} \end{picture}}}
\def \trivialknot{\raisebox{-4pt}{\begin{picture}(20,16)(-2,0)  
\thicklines
\put(8,8){\circle{16}} \end{picture}}}

\makeatother

\baselineskip=18pt

\title{\vspace*{-2cm}\large KNOT POLYNOMIALS AND GENERALIZED MUTATION}

\author{{\normalsize R.P.ANSTEE}\\
{\small Mathematics Department, University of British Columbia} \\
{\small 121-1984 Mathematics Road, Vancouver, BC, Canada, V6T1Y4}\\[-1mm]
\\[3mm]
{\normalsize J.H. PRZYTYCKI}\\
{\small Warsaw University, Warsaw, Poland}\\[-1mm]
\\[3mm]
{\normalsize D.ROLFSEN}\\
{\small University of British Columbia}\\
{\small 121-1984 Mathematics Road, Vancouver, 
BC, Canada, V6T1Y4}\\[-1mm]
}
\date{\empty}

\maketitle

\begin{abstract}
The motivation for this work was to construct a nontrivial knot with 
trivial Jones polynomial. Although that open problem has not yielded, 
the methods are useful for other problems in the theory of knot polynomials. 
The subject of the present paper is a generalization of Conway's mutation 
of knots and links. Instead of flipping a 2-strand tangle, 
one flips a many-string tangle to produce a generalized mutant. 
In the presence of rotational symmetry in that tangle, 
the result is called a ``rotant". We show  that if a rotant 
is sufficiently simple, then its Jones polynomial agrees 
with that of the original link. As an application, this provides 
a method of generating many examples of links with the same Jones 
polynomial, but different Alexander polynomials. 
Various other knot polynomials, as well as signature, 
are also invariant under such moves, if one imposes more 
stringent conditions upon the symmetries. 
Applications are also given to polynomials of satellites and symmetric knots.
\end{abstract}

{\small {AMS(MOS)Subj.Class.: 57M25}}\\
{\small{link polynomial, signature, Jones polynomial, mutant, knot, link, skein}}\\
\newpage
\section{Rotors and rotants} 
%%%%%%%%%section1.tex%%%rotor-and-rotants%%%%%%%%%%
~~~\par
The reader is assumed to be somewhat familiar with classical knot theory, standard references being \cite{APR3,APR5,APR23} but with \cite{APR11,APR13,APR9,APR17} the most relevant to what we are discussing. For the reader's convenience, the defining axioms for those knot polynomials which we consider are recalled in Section 2. The ideas of tangle and of mutation of knots and links were introduced by Conway \cite{APR4}. Mutation of a link is achieved by locating a tangle, with two inputs and two outputs, and flipping the tangle over like a pancake, or rotationg it $180^{\circ}$, as in Fig.1. Mutants can thus be formed in three ways; see \cite{APR17} for a precise definition.
\par\vs
\begin{center}
\begin{tabular}{cc} 
\includegraphics[trim=0mm 0mm 0mm 0mm, width=.25\linewidth]
{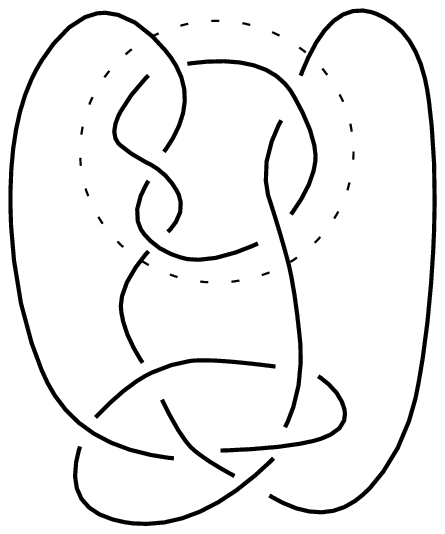}
&\hspace*{10mm}
\includegraphics[trim=0mm 0mm 0mm 0mm, width=.25\linewidth]
{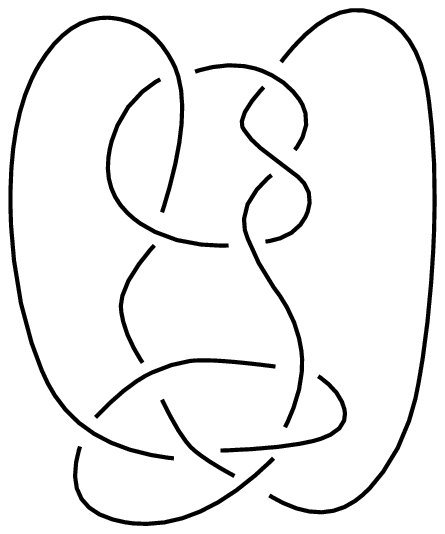}\\
\end{tabular}
\par\vs
Fig.1. Famous mutants: the Conway and Kinoshita-Terasaka knots.
\end{center}
\par\vs
The generalization of mutation which we consider here was inspired by the relation between the Jones polynomial of links and the Tutte 
dichromate polynomial for graphs \cite{APR12,APR13,APR24}, together 
with the ``rotor theorem" for the Tutte polynomial 
\cite{APR2,APR25,APR26}.
\par
Consider a planar projection $L$ of a classical link. If $B$ is a disk 
whose boundary is transverse to $L$, then $R=L \cap B$ is called a tangle. 
For our discussion, it will be convenient to picture $B$ as a metrically 
regular $n$-gon, $n \geq 3$, and $L$ is assumed to interesect the boundary 
of $B$ in exactly two points on each face, arranged symmetrically as 
in Fig.2. We call $R$ a {\it{rotor}} or {\it{$n$-rotor}}, if $R$ is 
(setwise) invariant under the rigid rotation $\rho: B \rightarrow B$ 
by angle $2 \pi $/$n$ around the center of $B$; the sense of overcrossings 
must also be respected. Following Tutte, we call the tangle $S$ 
complementary to $R$ in $L$ a {\it{stator;}} thus $L=S \cup R$. 
Now consider a line passing through the center of $B$ and either a corner 
or midpoint of an edge. Let $\mu:B \rightarrow B$ be a $\pi$-rotation 
(through the third dimension) with that line as axis. Thus $\rho$ and 
$\mu$ generate the dihedral group of symmetries of $B$. Although the 
boundary of the rotor enjoys full dihedral symmetry, we do not assume 
the rotor does $-$ indeed that would be an uninteresting case in what 
follows. Then $\mu R$ is a tangle with the same boundary as $R$. 
The crossings behave as if they were rigidly $3$-dimensional; that is 
a strand which passes over another will pass under, in the diagram, 
after being flipped by $\mu$. The generalized mutant $L'=S \cup \mu R$ 
will be called a {\it{rotant}} or {\it{$n$-rotant}} of the link $L$.
It is easy to see that $L$ and its rotant $L'$ have the same number 
of components (in fact, if two points of the boundary of $R$ are connected 
by a strand of $R$, then they are also connected by a strand of $\mu R$). 
Also, the choice of axis is irrelevant, and if one flips the stator 
instead of the rotor, the result is the same, up to ambient 
isotopy (even balanced isotopy, c.f. Lemma 2.1).
\par
\par\vs
\begin{center}
\begin{tabular}{c} 
\includegraphics[trim=0mm 0mm 0mm 0mm, width=.6\linewidth]
{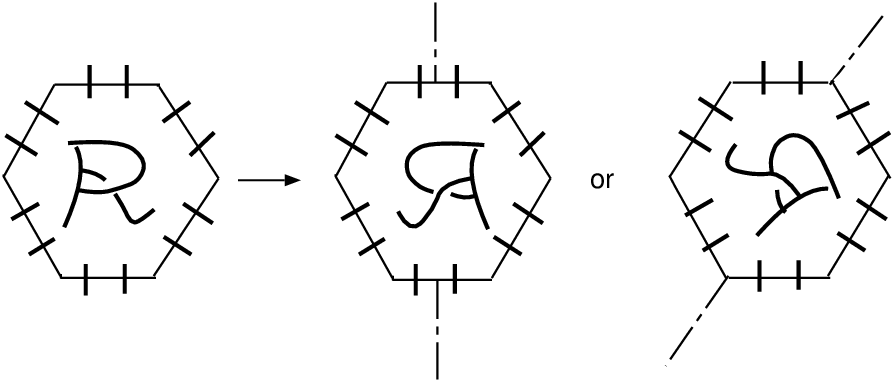}
\\
\end{tabular}
\par\vs
Fig.2. Allowable (and equivalent) flips to form a rotant.
\end{center}
\par\vs
So far we have made no assumptions about orientations, but if $L$ 
is considered as an {\it{oriented}} link we may require that $\rho$ 
either preserves all string orientations or reverses them all. 
It follows that $\mu$ either takes all inputs of the rotor to inputs 
or else carries all inputs to outputs. The {\it{oriented rotant L'}} 
is oriented by the original orientation of the stator, whereas all 
the strings of the rotor (besides getting flipped by $\mu$) keep or 
reverse their orientation according as $\mu$ preserves or reverses 
the sense of all in/outputs. We allow the possibility of closed 
strings in the rotor, and reverse their orientations as well in 
the latter case (see Fig. 3.).  
\par\vs
\begin{center}
\begin{tabular}{cc} 
\includegraphics[trim=0mm 0mm 0mm 0mm, width=.25\linewidth]
{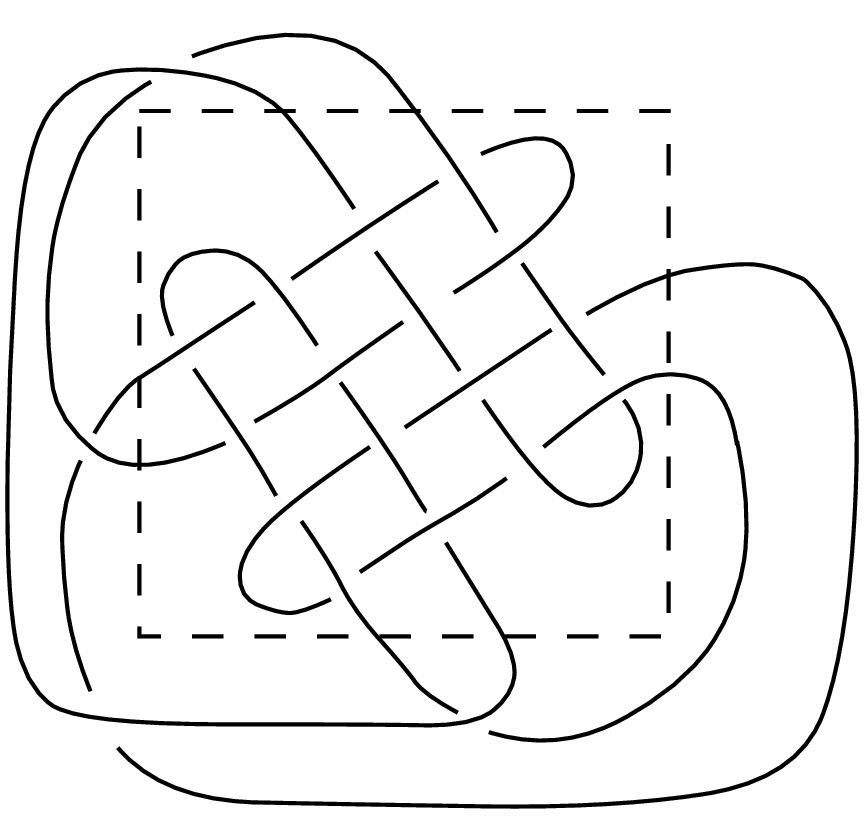}
&\hspace*{10mm}
\includegraphics[trim=0mm 0mm 0mm 0mm, width=.25\linewidth]
{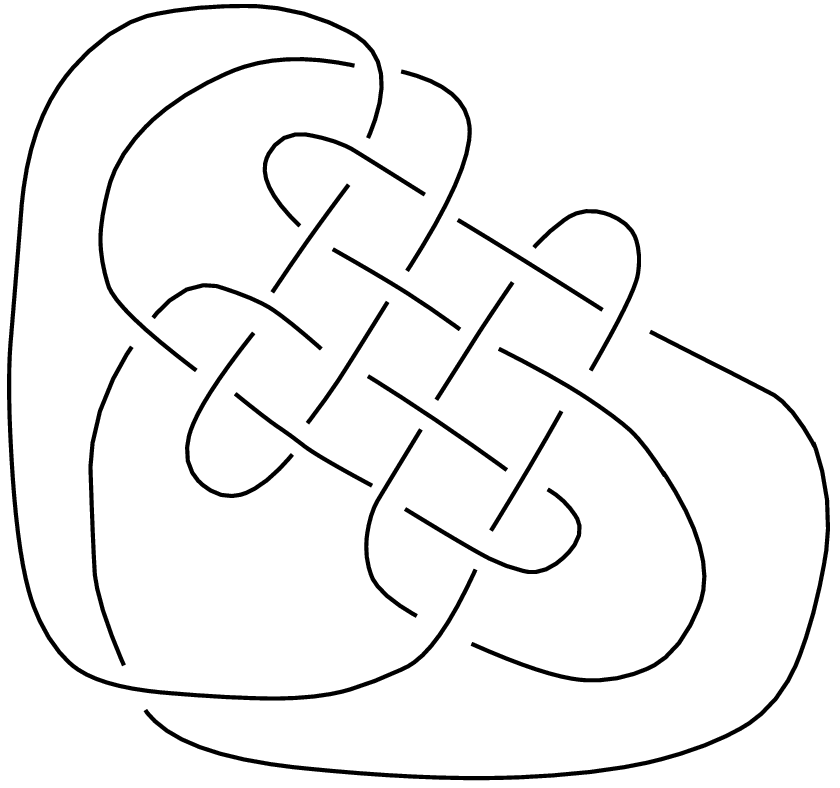}\\
\end{tabular}
\par\vs
Fig.3. A (non-oriented) rotor in link $L$, and the rotant $L'$.
\end{center}
\par\vs\noindent
{\bf{1.1. Example.}} As we can see from an example of Bleiler, 
Fig.4, certain conceivable flip must be disallowed if we are to prove 
any sort of polynomial invariance, namely those obtained by flipping 
on an axis passing through two boundary points. In this example, 
such a move transforms a trivial knot into a nontrivially linked set 
of three trefoils. Note that our definition of rotant rules this out. 
It also disallows mirror reflection of $R$.
\par\vs
\begin{center}
\begin{tabular}{cc}
\includegraphics[trim=0mm 0mm 0mm 0mm, width=.25\linewidth]
{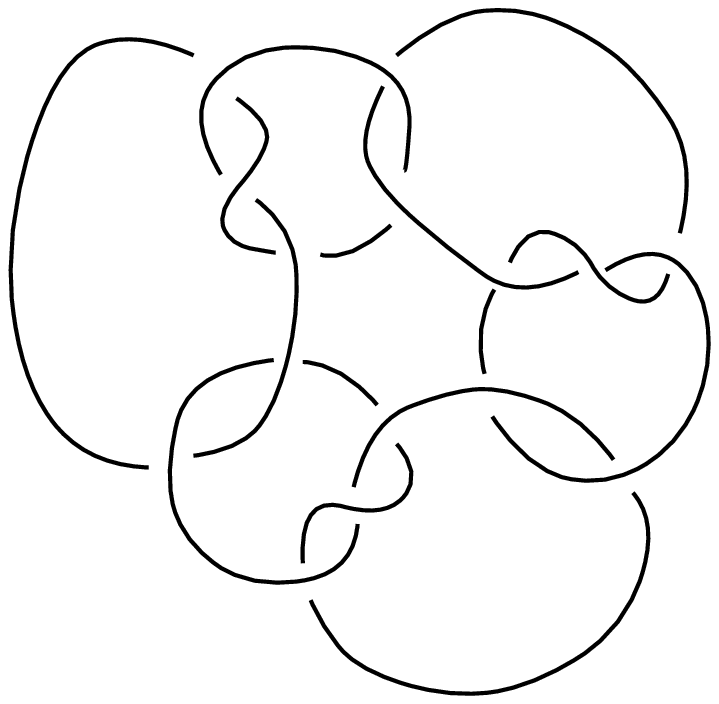}
&\hspace*{10mm}
\includegraphics[trim=0mm 0mm 0mm 0mm, width=.27\linewidth]
{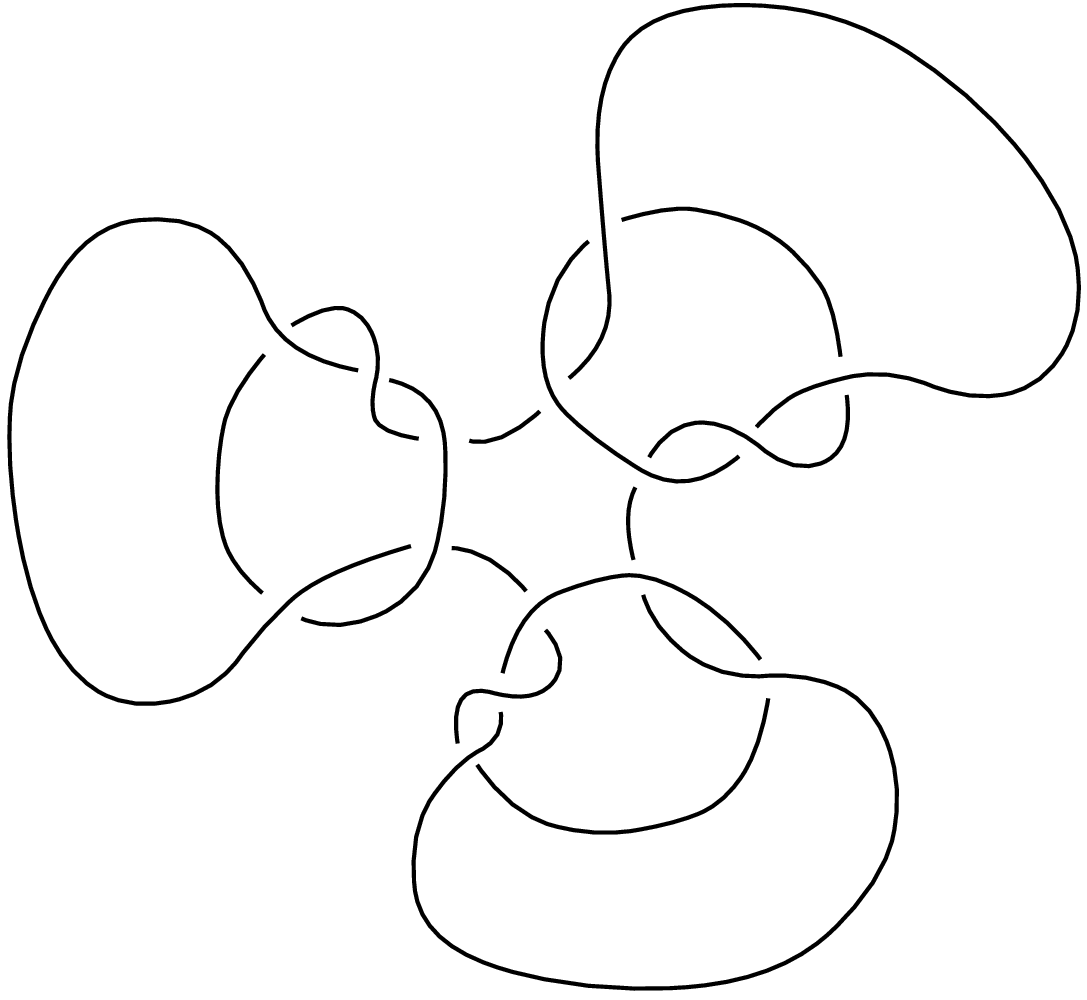}\\
\end{tabular}
\par\vs
Fig.4. These are not rotants of each other.
\end{center}

\par
\section{Link polynomials and resolving trees} 
%%%%%%%%%section2.tex%%%%%%%Link polynomials and resolving trees%%%%%%%%%%%%%%
~~~\par
To fix notation and define the notion of isomorphism of resolving 
trees we recall some of the recently discovered link-theoretic polynomials.
First we remind the reader that the relation of {\it{ambient isotopy}} 
(between link projections) is generated by the three Reidemeister moves 
$R1$, $R2$, $R3$ of Fig.5. The two moves $R2$ and $R3$ generate 
the more restrictive relation that Kauffman named {\it{regular isotopy}}. 
It will be natural for us to consider an intermediate relation, 
which we will call {\it{balanced isotopy}}, generated by moves $R2$, $R3$ 
and $BR1$, the ``balanced" move which introduces or deletes an opposite 
pair of curls, either in the same or in different components of a link.
\par\vs
\begin{center}
\begin{tabular}{cccc} 
$R1$ & $R2$ & $R3$ & $BR1$\\
\includegraphics[trim=0mm 0mm 0mm 0mm, width=.12\linewidth]
{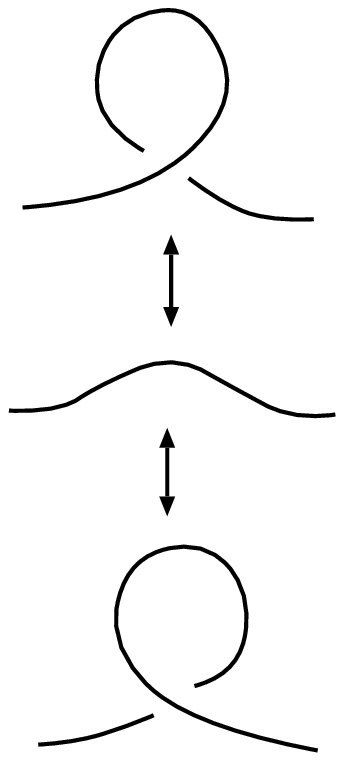}
&
\includegraphics[trim=0mm 0mm 0mm 0mm, width=.1\linewidth]
{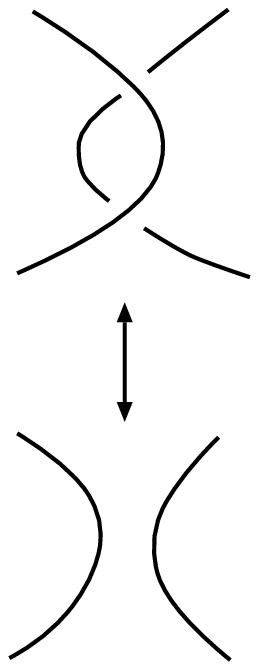}
&
\includegraphics[trim=0mm 0mm 0mm 0mm, width=.1\linewidth]
{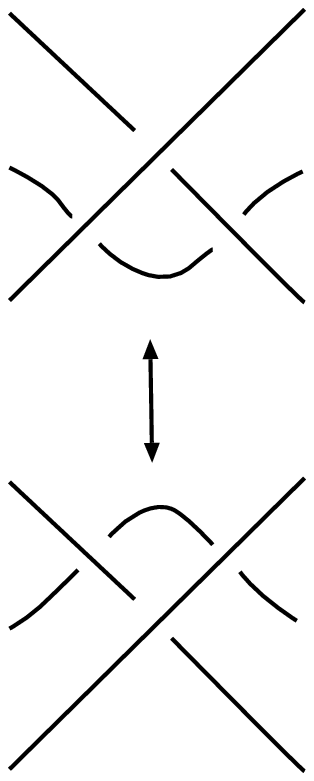}
&
\includegraphics[trim=0mm 0mm 0mm 0mm, width=.13\linewidth]
{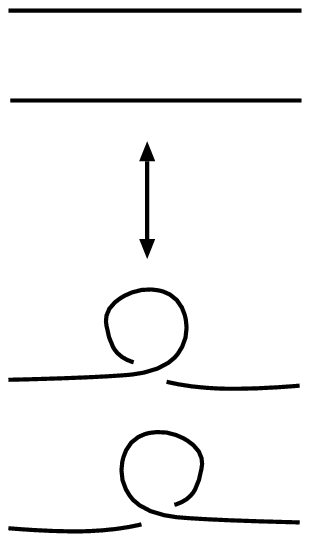}\\
\end{tabular}
\par\vs
Fig.5. Reidemeister moves.
\end{center}
\par\vs
It is not difficult to verify that two  oriented link diagrams $L$ and $L'$ are balanced isotopic if and only if they are ambient isotopic and w($L$)=w($L'$). Here w($L$) is the {\it{planar writhe}} defined by taking the algebraic sum of the crossings, counting 
%{\includegraphics[trim=0mm 0mm 0mm 0mm, width=.02\linewidth] {APRfigB-11.eps}}
{\includegraphics[trim=0mm 0mm 0mm 0mm, width=.03\linewidth] {L+maly.eps}}
~and~
%{\includegraphics[trim=0mm 0mm 0mm 0mm, width=.02\linewidth]{APRfigB-21.eps}} 
{\includegraphics[trim=0mm 0mm 0mm 0mm, width=.03\linewidth] {L-maly.eps}}
 as $+1$ and $-1$, respectively. 
For unoriented links $L$ and $L'$, the criterion for ambient isotopy to imply balanced isotopy is that sw($L)=$sw($L'$), where the {\it{self writhe}} sw($L$) is a similar sum, but counting only crossings involving two strings from the same component (this  is independent of how one orients $L$). The following lemma is a direct consequence of these facts.
\par\vs\noindent
{\bf{2.1.Lemma.}} {\it{Suppose L is a link diagram, corresponding}} ({\it{say}}) {\it{to the projection of a link in the \{x,y,z\}-space into the plane z=0. Apply the flip $(x,y,z) \rightarrow (-x,y,-z)$ and let $\mu L$ denote the projection of the resulting link. Then L and $\mu L$ are balanced isotopic.}}
\par\vs
Perhaps the simplest of the new polynomials is the Kauffman {\it{bracket}} 
$\langle L \rangle (A) \in \Bbb Z$[$A^{\pm}$] defined for unoriented 
link diagrams $L$, with defining relations, as in \cite{APR12}:
\par
$\langle {\includegraphics[trim=0mm 0mm 0mm 0mm, width=.02\linewidth]
{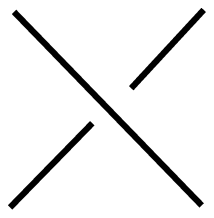}} \rangle =A \langle {\includegraphics[trim=0mm 0mm 0mm 0mm, width=.02\linewidth]
{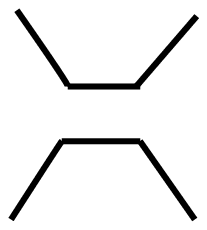}} \rangle +A^{-1} \langle {\includegraphics[trim=0mm 0mm 0mm 0mm, width=.02\linewidth]
{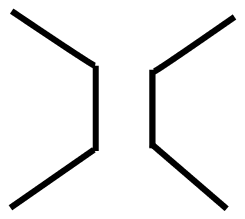}} \rangle,$ ~~~$\langle O \cup K \rangle =(-A^2-A^{-2}) \langle K \rangle,$ ~~~$\langle O \rangle =1.$   ~~~~~~~~~~(2.1)
\par\noindent
Here the symbols {\includegraphics[trim=0mm 0mm 0mm 0mm, width=.02\linewidth]
{APRfigA-4.eps}},~{\includegraphics[trim=0mm 0mm 0mm 0mm, width=.02\linewidth]
{APRfigA-2.eps}}~ 
and ~{\includegraphics[trim=0mm 0mm 0mm 0mm, width=.02\linewidth]
{APRfigA-3.eps} stand for links which look like that in a neighborhood 
of a point and are identical elsewhere, $O$ is a curve in the diagram 
with no crossing points and $\cup$ is disjoint union. The bracket 
polynomial is an invariant of regular isotopy and, in fact, 
of balanced isotopy.
\par
Calculation of $\langle L \rangle $ for an unoriented link diagram 
may be accomplished inductively, using the first equation of (2.1) 
at a crossing to express $\langle L \rangle $ in terms of 
the bracket polynomial of two simpler link diagrams. 
The calculation may be recorded by a binary tree, which we will call 
a {\it{bracket resolving tree}}. Each node of the tree corresponds 
to a link diagram, and if it has successors they stand in the relation:
\par\vs
\begin{center}
\begin{tabular}{c} 
\includegraphics[trim=0mm 0mm 0mm 0mm, width=.35\linewidth]
{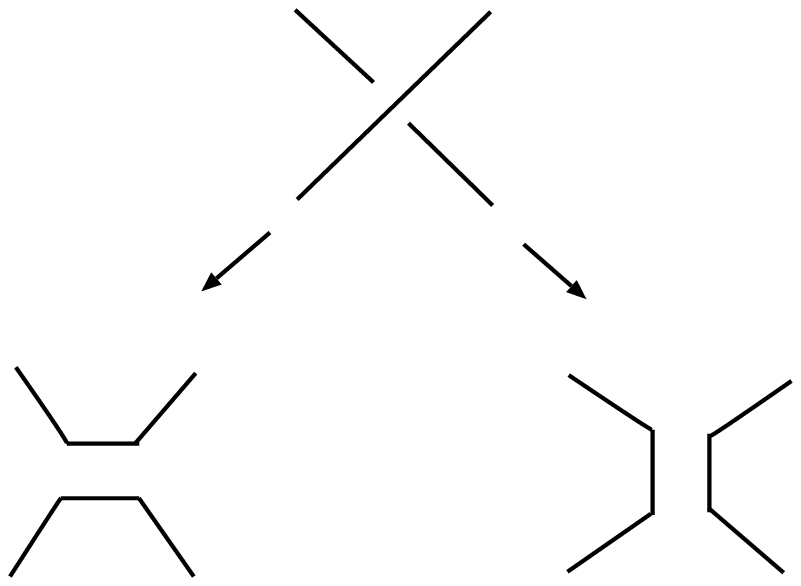}
\end{tabular}
\end{center}
\par\vs\noindent
The terminal nodes are the {\it{leaves}} of the tree. 
One may always build a bracket resolving tree whose leaves have 
no crossings, but we do not require this.
Two bracket resolving trees will be called {\it{isomorphic}} if 
\par
(1) they are isomorphic as binary trees,
\par
(2) for each non-terminal node which, with its two successors stand 
in the relation {\includegraphics[trim=0mm 0mm 0mm 0mm, width=.02\linewidth]
{APRfigA-4.eps}},~ {\includegraphics[trim=0mm 0mm 0mm 0mm, width=.02\linewidth]
{APRfigA-2.eps}~ and~ 
{\includegraphics[trim=0mm 0mm 0mm 0mm,
width=.02\linewidth]{APRfigA-3.eps}, the corresponding nodes under the
isomorphism stand in the same relation, and 
\par
(3) the corresponding leaves under the isomorphism represent link diagrams which are {\it{balanced}} isotopic.
\par
One can then define an  invariant of unoriented links under ambient isotopy by the formula
\par
${\underbar{\it{f}}}_{L}(A)=(-A)^{-3{\rm{sw}}(L)} \langle L \rangle$.
\par\noindent
Alternatively one can form an invarinant of oriented links under 
ambient isotopy: orient the link and define
\par
$f_L(A)=(-A)^{-3{\rm{w}}(L)} \langle L \rangle $.
\par\noindent
This is a version of the Jones polynomial \cite{APR11,APR13}: 
$V_L(t)=f_L(t^{-1/4})$.
\par
The several-variable {\it{skein}} polynomial (also called HOMFLY,
FLYPMOTH, twisted Alexander, Jones-Conway, generalized Jones, etc.) has
several equivalent versions, as in \cite{APR6,APR22}. We use Hoste's
notation \cite{APR9}. A (Laurent) polynomial $P_L(x,y,z)$, for each
$oriented$ link $L$ is defined by the equations :
\par
$xP_{{\includegraphics[trim=0mm 0mm 0mm 0mm,
width=.02\linewidth]{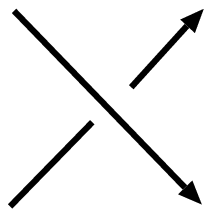}}}+
yP_{{\includegraphics[trim=0mm 0mm 0mm 0mm,
width=.02\linewidth]{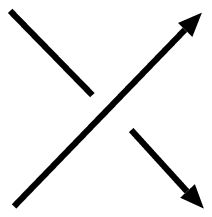}}}+
zP_{{\includegraphics[trim=0mm 0mm 0mm 0mm,
width=.02\linewidth]{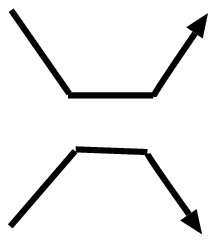}}}=0$ and $P_{\rm{unknot}}=1$. ~~~~~~~~~~~~~~~~~~~~~~~~~~~~~~~~~   (2.2)
\par
The Lickorish-Millett polynomial $P_L(l,m)$ \cite{APR17} is related 
by the substitution $x=l, y=l^{-1},z=m$. As above, one can 
compute $P_L$ for an oriented diagram by a binary calculation 
tree which we will call a {\it{skein resolving tree}}. 
Each non-terminal node, with its two successors, stands in the relation
\par\vs
\begin{center}
\begin{tabular}{c} 
\includegraphics[trim=0mm 0mm 0mm 0mm, width=.35\linewidth]
{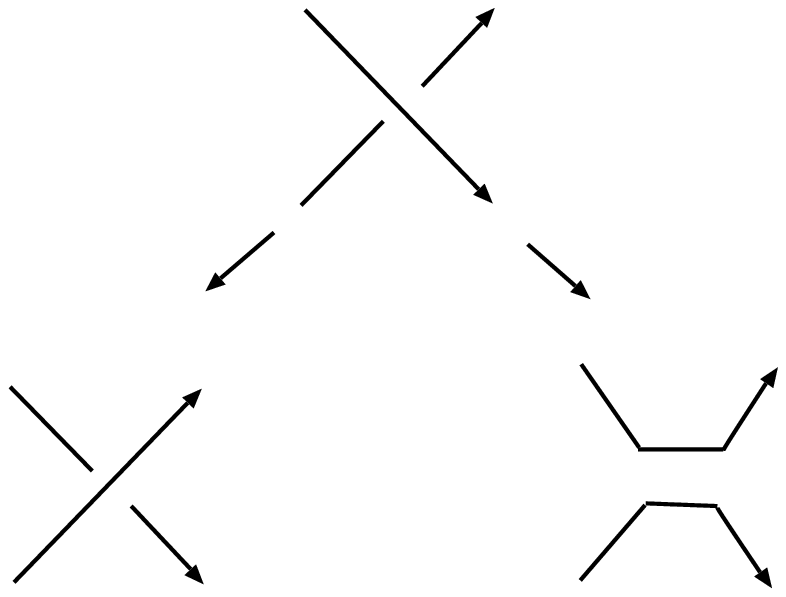}\\
~~\\
or else \\
~~\\
\includegraphics[trim=0mm 0mm 0mm 0mm, width=.35\linewidth]
{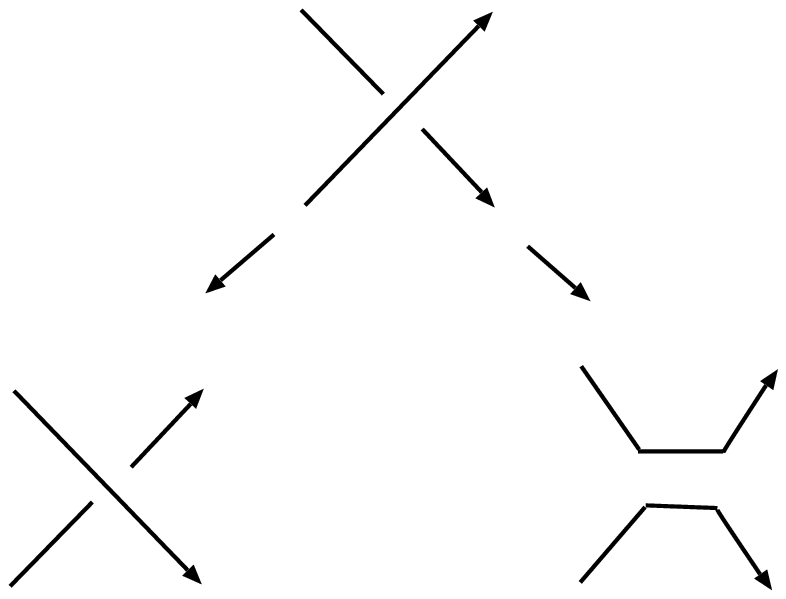}\\
\end{tabular}
\end{center}
\par\vs\noindent
and a skein tree isomorphism is assumed to preserve such relations.
\par
Finally, the leaves, under an isomorphism, are assumed to be ambient 
isotopic as oriented link diagrams. Skein equivalence, 
as exposited in \cite{APR17}, is the equivalence relation of oriented 
links generated by isomorphism of skein trees. 
It is an open question whether skein equivalent links actually 
have isomorphic skein resolving trees.
\par
Now we turn to the {\it{Kauffman}} 2-variable polynomial, 
${\Lambda}_L(a,x)$ \cite{APR13}, defined for unoriented link diagrams 
by the axioms:
\par
${\Lambda}_{{\includegraphics[trim=0mm 0mm 0mm 0mm, width=.02\linewidth]
{APRfigA-4.eps}}} +{\Lambda}_{{\includegraphics[trim=0mm 0mm 0mm 0mm, 
width=.02\linewidth]
{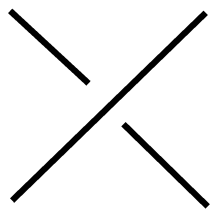}}}=
x{\Lambda}_{{\includegraphics[trim=0mm 0mm 0mm 0mm, width=.02\linewidth]
{APRfigA-2.eps}}}+
x{\Lambda}_{{\includegraphics[trim=0mm 0mm 0mm 0mm, width=.02\linewidth]
{APRfigA-3.eps}}}$,
\par
${\Lambda}_{{\includegraphics[trim=0mm 0mm 0mm 0mm, width=.02\linewidth]
{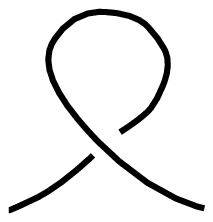}}}=
a{\Lambda}_{{\includegraphics[trim=0mm 0mm 0mm 0mm, width=.02\linewidth]
{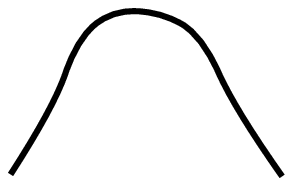}}}$,~ ${\Lambda}_{{\includegraphics[trim=0mm 0mm 0mm 0mm, width=.02\linewidth]
{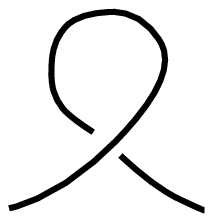}}}=
a^{-1}{\Lambda}_{{\includegraphics[trim=0mm 0mm 0mm 0mm, width=.02\linewidth]
{APRfigarc.eps}}}$,~${\Lambda}_O=1$.~~~~~~~~~~~~~~~~~~~~~~~~~~~~~~~~~~~~~    (2.3)
\par\noindent
${\Lambda}$ is invariant under balanced isotopy.
\par
A {\it{Kauffman resolving tree}} is a ternary tree such that each node 
which is not a leaf stands, with its three successors, in the relation:
\par\vs
\begin{center}
\begin{tabular}{c} 
\includegraphics[trim=0mm 0mm 0mm 0mm, width=.35\linewidth]
{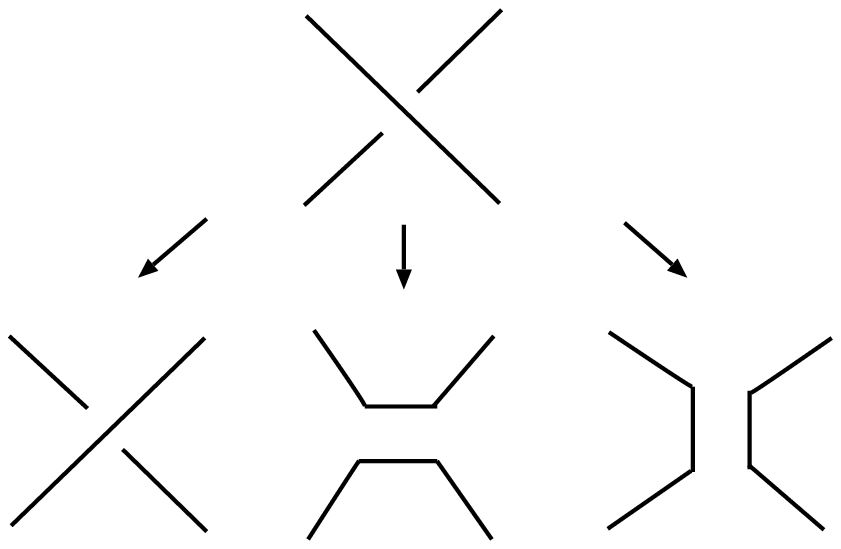}
\end{tabular}
\end{center}
\par\vs\noindent
An isomorphism of Kauffman trees must respect this relation. 
(Note that confusing the sequence 
{\includegraphics[trim=0mm 0mm 0mm 0mm, width=.02\linewidth]
{APRfigA-4.eps}},~{\includegraphics[trim=0mm 0mm 0mm 0mm, 
width=.02\linewidth]
{APRfigA-1.eps}},~{\includegraphics[trim=0mm 0mm 0mm 0mm, 
width=.02\linewidth]
{APRfigA-2.eps}},~{\includegraphics[trim=0mm 0mm 0mm 0mm, width=.02\linewidth]
{APRfigA-3.eps}}~ 
with~ {\includegraphics[trim=0mm 0mm 0mm 0mm, width=.02\linewidth]
{APRfigA-4.eps}},~{\includegraphics[trim=0mm 0mm 0mm 0mm, width=.02\linewidth]
{APRfigA-1.eps}},~{\includegraphics[trim=0mm 0mm 0mm 0mm, width=.02\linewidth]
{APRfigA-3.eps}},~{\includegraphics[trim=0mm 0mm 0mm 0mm, width=.02\linewidth]
{APRfigA-2.eps}} would not affect the calculation of $\Lambda$. 
Nevertheless, we do not permit this in the definition since our 
results actually ensure the stronger notion of isomorphism.) 
The diagrams associated with corresponding leaves, under an isomorphism, 
are assumed to be related by {\it{balanced}} isotopy.
\par
The Kauffman invariant of ambient isotopy is defined by giving 
$L$ an orientation and defining 
\par
~~~$F_L(a,x)=a^{-{\rm{w}}(L)}{\Lambda}_L(a,x)$.
\par
As above one can define a similar invariant for unoriented links by 
\par
~~~${\underbar{\it{F}}}_L(a,x)=a^{-{\rm{sw}}(L)}{\Lambda}_L(a,x)$.
\par\noindent
(Our results regarding $f$ and $F$ have corresponding versions for
${\underbar{{\it{f}}}}$ and ${\underbar{{\it{F}}}}$ which we omit in the
interest of simplicity.) The predecessor of Kauffman's polynomial is the
polynomial of Brandt-Lickorish-Millett \cite{APR1} and Ho \cite{APR8},
which is just $Q_L(x)=F_L(1,x)={\Lambda}_L(1,x)$. It does not depend on
orientation. Finally, we remind the reader that the Conway and Alexander
polynomials are, respectively, ${\bigtriangledown}_L(z)=P_L(1,-1,z)$ and
${\Delta}_L(t)={\bigtriangledown}_L(t^{1/2}-t^{-1/2})$, and the Jones
polynimial satisfies
$V_L(t)=P_L(t^{-1},-t,t^{-1/2}-t^{1/2})=F_L(t^{-3/4},-t^{1/4}-t^{-1/4})$.
See [14 or 13] for details.

\par
\section{Main results}
%%%%%%%%%%%%%%%%3.Main results%%%section3.tex%%%%%%%%%%%%%%%
~~~
\par\noindent
{\bf{3.1. Theorem.}} {\it{Suppose $n \leq 5$. If the link $L'$ is an
$n$-rotant of $L$, then they have isomorphic bracket resolving trees and
hence have the same Kauffman bracket}} ({\it{and, in the case of knots,
equal Murasugi signatures}}). {\it{If they are oriented $n$-rotants, 
they also have the same Jones polynomial.}}
\par\vspace{1cm}\noindent
{\bf{3.2. Theorem.}} {\it{If $n \leq 4$ and $L'$ is an oriented $n$-rotant 
of $L$, then they have isomorphic skein resolving trees, 
so $P_L(x,y,z)=P_{L'}(x,y,z)$ and their Alexander, Conway and 
Jones polynomials also agree.}}
\par\vspace{1cm}\noindent
{\bf{3.3. Theorem.}} {\it{If $L'$ is a $3$-rotant of $L$, then they 
have isomorphic Kauffman resolving trees and 
${\Lambda}_L(a,x)={\Lambda}_{L'}(a,x)$. If they are oriented $3$-rotants, 
then also they have the same Kauffman polynomials}}: {\it{$F_L(a,x)=
F_{L'}(a,x)$}}.
\par\vspace{1cm}
Before proving the theorems, we should point out possible 
application in connection with questions which are (at this writing) 
still open.
\par\vspace{1cm}\noindent
{\bf{3.4. Question.}} Is there a nontrivial knot with trivial 
Jones polynomial?
More generally, is there a nontrivial link whose Jones polynomial 
equals that of an unlink?\footnote{Added for e-print: The question for 
knots remains open. For links
there are examples of nontrivial links whose Jones polynomial equals
that of an unlink: M.B.Thistlethwaite,
Links with trivial Jones polynomial,
{\it J. Knot Theory and its Ramifications} 10(4), 2001, 641-643.\ 
S.Eliahou, L.H.Kauffman and M.B.Thistlethwaite,
Infinite families of links with trivial Jones polynomial,
Topology, 42 (2003) 155-169.}
\par\vspace{1cm}
An affirmative answer would be implied by one for the following, for $n \leq 5$.\par\vspace{1cm}\noindent
{\bf{3.5. Question.}} Does the unknot (unlink) have an $n$-rotant which is knotted (non-trivially linked)?
\par\vspace{1cm}
In fact, if there is a nontrivial $3$-rotant of an unknot, there are infinitely many different prime knots with trivial Jones polynomial (see Remark 4.6). Similarly, one might hope to find, via a rotant of an unknot, a method of finding nontrivial knots with trivial $P,Q$ or $F$ polynomial.
\par\vspace{1cm}\noindent
{\bf{Proof of Theorem 3.1.}} Consider bracket resolving trees for $L$ and $L'$, in which one eliminates only crossings of the stator, without altering the rotor tangle. Clearly the trees are isomorphic except possibly for the leaves, which correspond to links in which the stator has no crossings. For $n=5$, there are just the ten possibilities of Fig.6, up to rotation. Each has an axis of symmetry, and so we see that, in the leaves, forming the rotant changes corresponding links by a rigid $180^{o}$ rotation through space. But such links are balanced isotopic, by Lemma 2.1. The cases $n \leq 4$ are similar and left to the reader. The second part of Theorem 3.1 follows, since the writhe is also unchanged by forming an oriented rotant. Finally, \cite[Theorem 4.5]{APR21} implies that if two knots have isomorphic bracket resolving trees, then they have the same signature.
\par\hfill$\Box$
\par\vs
\begin{center}
\begin{tabular}{c} 
\includegraphics[trim=0mm 0mm 0mm 0mm, width=.7\linewidth]
{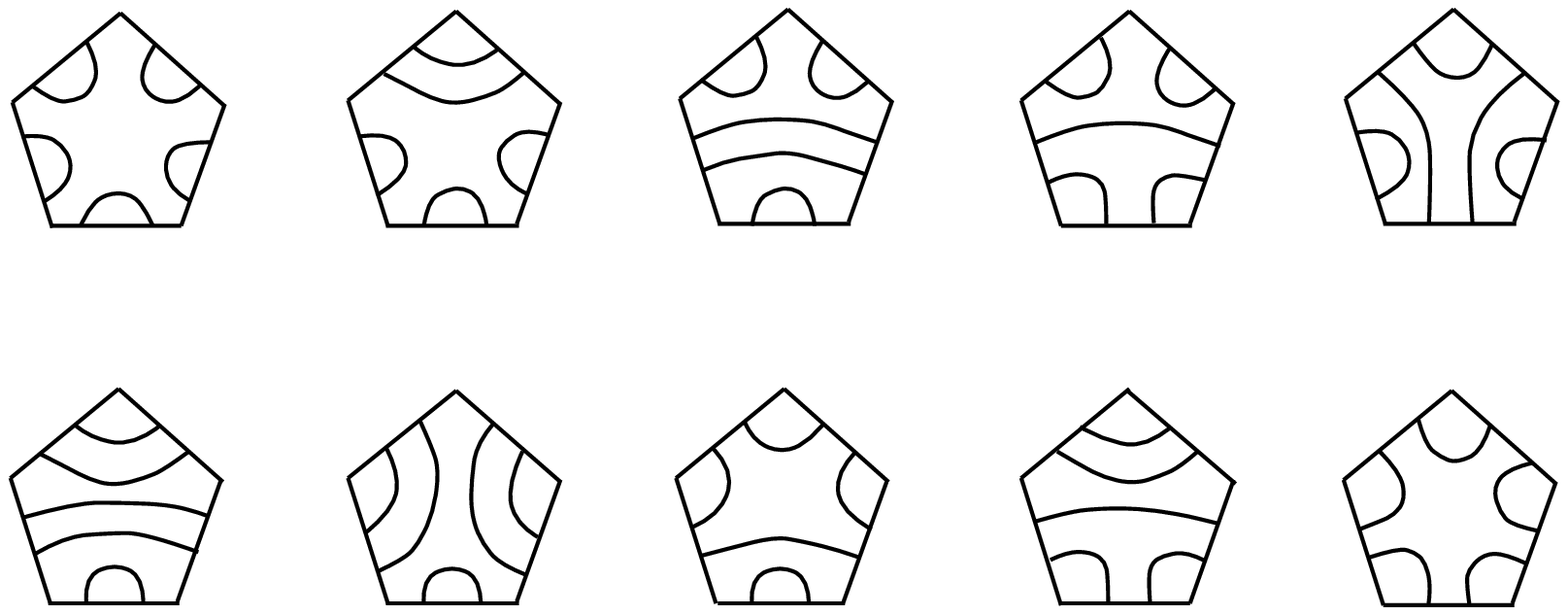}\\
\end{tabular}
\par\vs
Fig.6. The possible leaf stators in Theorem 3.1.
\end{center}
\par\vs
\par\vspace{1cm}\noindent
{\bf{Proof of Theorem 3.2 and 3.3.}} 
One may argue analogously to the proof of Theorem 3.1, but for variety 
we will prove Theorem 3.2 using the point of view of skein theory, 
as elaborated in \cite{APR17}. Consider $n=4$. For a fixed $4$-rotor $R$, 
the function $S \rightarrow (R \cup S)-(\mu R \cup S)$ defines a 
linear mapping ${\cal L}({\cal R}_4) \rightarrow {\cal L}(S^3)$ from 
the linear skein module of the room ${\cal R}_4$ of Fig.7 into the 
linear skein module of $S^3$. Under any choice of $4$ inputs and 
$4$ outputs, ${\cal L}({\cal R}_4)$ is generated by $4!=24$ ``standard" 
tangles and all have an axis of symmetry, except a few which, up to 
dihedral symmetry, are of five types shown in Fig.7. However, it is easy 
to check that none of these can occur as stator with an oriented rotor. 
Thus, the linear function vanishes on a set of generators for 
the linear skein of stators $S$ compatible with an oriented $R$, 
hence is identically zero there, proving Theorem 3.2. 
To prove Theorem 3.3 one only need to note that all the generators 
of the ``unoriented skein" theory of ${\cal R}_3$ may be chosen 
with axis of symmetry.
\hfill$\Box$
\par\vs
\begin{center}
\begin{tabular}{cc} 
\includegraphics[trim=0mm 0mm 0mm 0mm, width=.25\linewidth]
{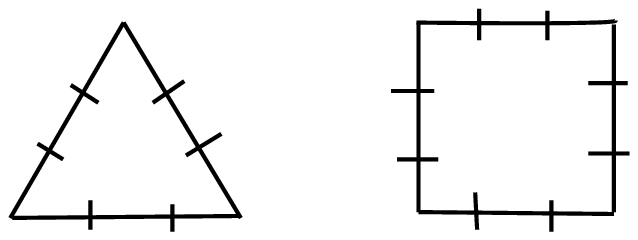}
&
\includegraphics[trim=0mm 0mm 0mm 0mm, width=.7\linewidth]
{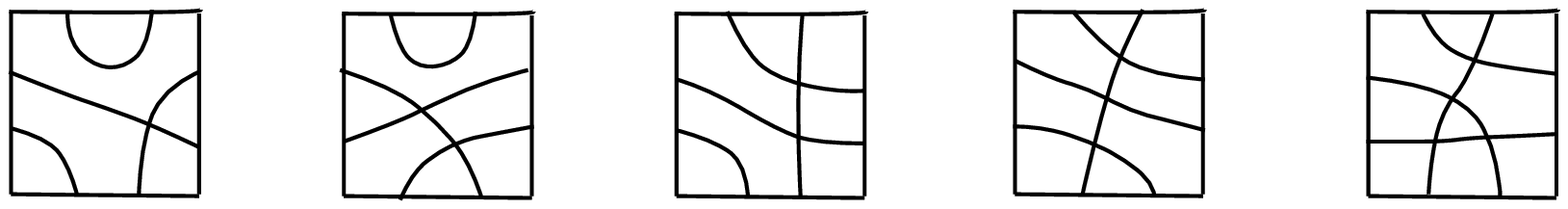}
\\
\end{tabular}
\par\vs
Fig.7. Rooms ${\cal R}_3,{\cal R}_4$ and asymmetric generators for ${\cal L}({\cal R}_4)$.
\end{center}
\par\vs\noindent
{\bf{3.6. Example.}} The 4-rotants $L$ and $L'$ of Fig. 3 have 
the same Jones polynomials: 
\par
$V_L(t)=V_{L'}(t)=-6t^{-7}+29t^{-6}-84t^{-5}+178t^{-4}-298t^{-3}+422t^{-2}-514t^{-1}+550-522t+435t^2-314t^3+192t^4-94t^5+34t^6-8t^7+t^8$.
\par\vs
Their skein polynomials are quite different. Note that they are not oriented rotants, else by Theorem 3.2 the skein polynomials would agree. Even their Alexander and Conway polynomials disagree, the latter being
\par
${\bigtriangledown}_L(z)=1-2z^4-10z^6-8z^8-5z^{10}$,
\par\noindent
and
\par
${\bigtriangledown}_{L'}(z)=1-10z^4-8z^6+5z^8+6z^{10}+2z^{12}$.

\par
\section{Satellites, double rotors and further results}
%%%%%section4.tex%%%Satellites,double rotors and further results%%%
~~~\par
For certain types of rotors, one can improve upon the above results.
An $n$-rotor $R$ is called {\it{cup-trivial}} if, after replacing an 
adjacent pair of boundary points by a ``cup", the resulting tangle 
is isotopic (rel boundary) to a trivial tangle, as shown in Fig.8. 
\par\noindent
(Additional trivial components are allowed.)
\par\vs
\begin{center}
\begin{tabular}{cc} 
\includegraphics[trim=0mm 0mm 0mm 0mm, width=.65\linewidth]
{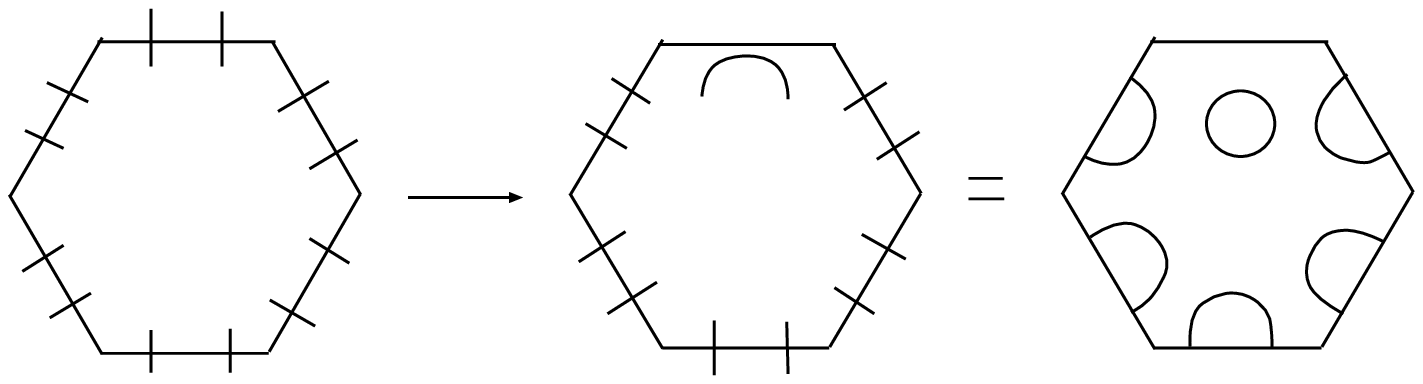}
&
\includegraphics[trim=0mm 0mm 0mm 0mm, width=.17\linewidth]
{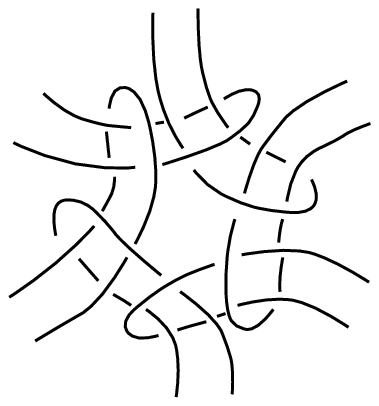}
\\
\end{tabular}
\par\vs
Fig.8. Cup-trivial rotor, and an example.
\end{center}
\par\vs\noindent
{\bf{4.1. Theorem.}} {\it{Suppose $L$ is an unoriented link diagram with a cup-trivial $n$-rotor, and let $L'$ be its corresponding rotant. Then}}
\par
({\it{a}}) {\it{for $n \leq 7$, $L$ and $L'$ have isomorphic bracket resolving trees and $\langle L \rangle = \langle L' \rangle$, and}} 
\par
({\it{b}}) {\it{if the rotor is oriented and $n \leq 5$, $L$ and $L'$
have isomorphic skein trees and}}
\par
{\it{consequently $P_L=P_{L'}$}}.
\par\vspace{1cm}\noindent
{\bf{Proof.}} Construct, as usual, resolving trees for $L$ and $L'$ 
changing only the common stator part. Then, in the leaves of the trees 
we see that each stator either has a ``cup" or an axis of symmetry. 
(The rather lengthy case-checking is left to the diligent reader.) 
In either case we conclude that they have isomorphic computation trees.
\hfill$\Box$
\par\vspace{1cm}
A link diagram $L$ which contains an $n$-rotor $R$ such that the stator 
$S$ is also an $n$-rotor is called a {\it{double $n$-rotor}}. 
A {\it{generalized n-rotor of type}} $k$ is defined in the same way 
as an $n$-rotor, except that one has $2k$ boundary points on each 
of the $n$ lateral faces, as in Fig.9. If two generalized $n$-rotors 
(of the same type) are put together, the resulting link diagram 
will be called a {\it{generalized double $n$-rotor}}.
\newpage
\begin{center}
\begin{tabular}{cc} 
\includegraphics[trim=0mm 0mm 0mm 0mm, width=.2\linewidth]
{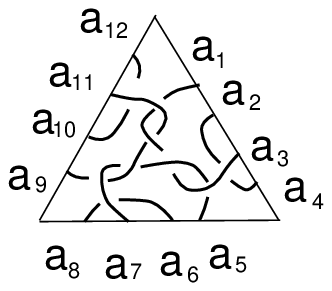}
&
\includegraphics[trim=0mm 0mm 0mm 0mm, width=.2\linewidth]
{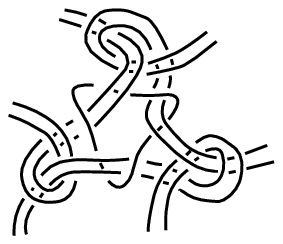}
\\
\end{tabular}
\par\vs
Fig.9. A generalized 3-rotor of type 2 and a parallel one.
\end{center}
\par\vs\noindent 
{\bf{4.2.Theorem.}} {\it{Let $n$ be an arbitrary positive integer and let $n^*=p_1 \cdots p_k$, where $n={p_1}^{a_1} \cdots {p_k}^{a_k}$ is written as product of powers of distinct primes.}}
\par
({\it{a}}) {\it{Suppose $L$ is a generalized double $n$-rotor and $L'$
is its rotant. Then if one}} 
\par
{\it{reduces coefficients modulo $n^*$, their Kauffman brackets agree}}:
\par
$\langle L \rangle (A) \equiv \langle L' \rangle (A) ($mod $n^*)$.}}
\par~~\par
({\it{b}}) {\it{If $L$ is a double $n$-rotor and $L'$ its rotant, then}}
\par
${\Lambda}_L(a,x) \equiv {\Lambda}_{L'}(a,x) ($mod $n^*),$
\par
{\it{and consequently}}
\par
$Q_L(x) \equiv Q_{L'}(x) ($mod $n^*).$
\par~~\par
({\it{c}}) {\it{If $L$ is an oriented double $n$-rotor and $L'$ its oriented rotant, then}}
\par
$P_L(x,y,z) \equiv P_{L'}(x,y,z) ($mod $n^*)$,
\par
{\it{and}}
\par
$F_L(a,x) \equiv F_{L'}(a,x)$ (mod $n^*)$,
\par
{\it{and therefore their Alexander, Conway and Jones polynomials also
agree,}}
\par
{\it{modulo $n^*$.}}
\par\vspace{1cm}\noindent
{\bf{Sketch of the proof.}} It is enough to prove the theorem for $n=p$, 
a prime. Following an observation of Murasugi \cite{APR19}, in a 
calculation tree any asymmetric link which appears, actually 
appears $p$ times and so the total contribution of such terms disappears 
when considering coefficients modulo $p$. Then one constructs a ``symmetric" 
calculation tree by altering crossings of the stator (which is itself a rotor) 
equivariantly until in each of the leaves the stator has axis symmetry.
\par\vs
We now turn to the question of satellites, to obtain results for 
cables of certain rotants, analogous to those of \cite{APR16,APR18,APR21} 
for cables of mutants.
\par\vs\noindent
{\bf{4.3. Definition.}} Consider a generalized $3$-rotor of type $2$ and 
label the boundary points consecutively $a_1,\cdots,a_{12}$ as in Fig.9. 
Call it {\it{parallel}} if each pair $a_{2i-1},a_{2i}$ is connected 
to a pair $a_{2j-1},a_{2j}$ by parallel strands. More precisely, 
these strands are assumed to be the boundary of three disjoint bands 
in the diagram. Moreover, we assume the bands are untwisted, in the 
sense of being ``flat", i.e. immersed in the projection plane. (As twists 
can be pushed into the stator, we are really assuming they are 
equally twisted.) Other closed curves are allowed, as long as they 
do not intersect the bands.
\par\vs\noindent
{\bf{4.4. Theorem.}} {\it{Let $L$ be an unoriented link diagram composed 
of a parallel generalized $3$-rotor of type $2$}} ({\it{as described 
above}}) {\it{and an arbitrary stator $S$, and let $L'$ be the corresponding 
rotant of $L$. Then $L$ and $L'$ have isomorphic bracket resolving 
trees and $\langle L \rangle =\langle L' \rangle$}}.
\par\vs\noindent
{\bf{Proof.}} As usual, build the bracket resolving trees 
for $L$ and $L'$ by simultaneously smoothing crossings only in 
the common stator $S$. Now consider a stator part of one of the leaves. 
If no string of this stator joins a point $a_{2i-1}$ with $a_{2i}$, 
then one can easily check that this stator has an axis of symmetry 
and the corresponding leaves of the two trees are equivalent. 
Now consider the case that $a_{2i-1}$ is joined with $a_{2i}$ in a leaf 
stator (by a cup). Then, using the fact that the rotor $R$ is parallel, 
we can simplify the stator essentially by a regular isotopy of the link 
in that leaf. However this simplification destroys the rotational 
symmetry of the rotor part. Observe, nevertheless, that if a pair 
of boundary points is joined to another pair of boundary points 
by parallel strings in the rotor, then their images under any 
axis symmetry are likewise joined by parallel strings in the rotor. 
If after simplifications of the leaf stator, one gets a stator 
with axis symmetry we are done (the simplified rotor does not 
depend on whether we used the leaf stator to simplify it or its image 
under an axis symmetry). We show that this is the case in Fig.10, under 
assumption that $a_1$ is joined to $a_2$ in the leaf stator 
and $a_1,a_2$ are joined by parallel strings in the rotor 
to $a_{11},a_{12}$. Other cases are left to the reader.
\hfill$\Box$
\par\vs\noindent
\par\vs
\begin{center}
\begin{tabular}{c} 
\includegraphics[trim=0mm 0mm 0mm 0mm, width=.8\linewidth]
{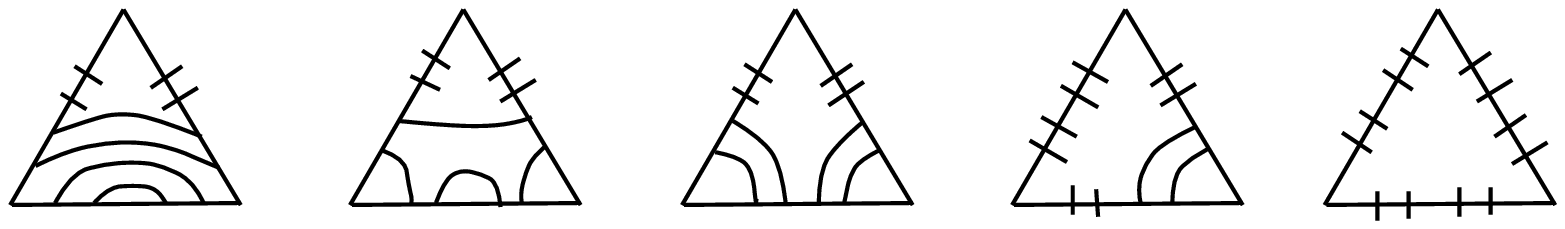}
\\
\end{tabular}
\par\vs
Fig.10. Simplified leaf stators.
\end{center}
\par\vs\noindent
{\bf{4.5.Corollary.}} {\it{Suppose an unoriented link diagram $K$ has $3$-rotant $K'$. Suppose further that $L$ and $L'$ are satellites, of wrapping number $2$}}  ({\it{e.g. $(2,q)$-cables or Whitehead doubles}}) {\it{of $K$ and $K'$, respectively}} ({\it{constructed in the same way}}). {\it{Then $L$ and $L'$ have isomorphic bracket resolving trees. In particular they have the same Kauffman bracket and their Jones polynomials differ only by a factor $t^{3r}$}}.
\par\vs\noindent
{\bf{Proof.}} By ``constructed in the same way" we mean that the satellite 
rotor forms parallel strands in the sense of Definition 4.3 and 
that the parts of $L$ and $L'$ in the stator are identical. 
Then it follows immediately.
\hfill$\Box$
\par\vs\noindent
{\bf{4.6. Remark.}} This shows that if one could find a $3$-rotant $K$ of 
the unknot, with $K$ knotted, then one could construct infinitely many 
distinct prime knots with trivial Jones polynomial by taking 
the various doublings as shown in Fig.11.
\par\vs
\begin{center}
\begin{tabular}{ccc} 
\includegraphics[trim=0mm 0mm 0mm 0mm, width=.5\linewidth]
{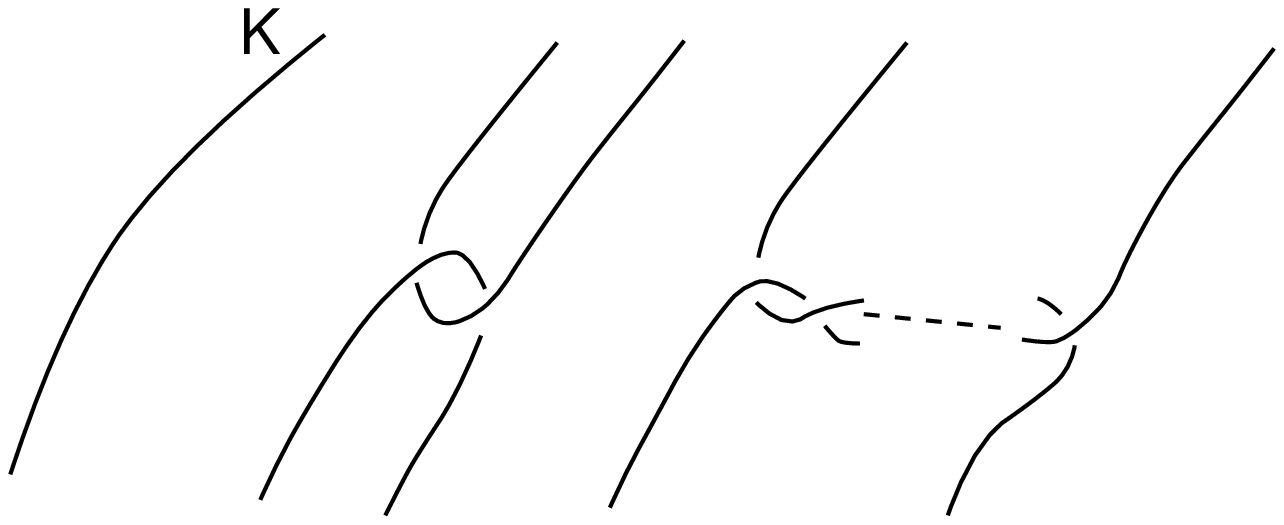}
\\
\end{tabular}
\par\vs
Fig.11. Various doublings of $K$.
\end{center}
\par\vs
Several of our results concerning the skein polynomials of rotants 
can be sharpened if we impose an additional condition on the stator. 
This was inspired by recent work of Jaeger \cite{APR10}.
\par\vs\noindent
{\bf{4.7.Definition.}}  Let $D$ be a part of a diagram of an oriented 
link (for example, a tangle or all of the link). 
We say that $D$ is a {\it{matched}} diagram if one can pair up the 
crossings in $D$ so that each pair is connected in the diagram 
as in Fig.12. Examples are given in Fig.13.
\par\vs
\begin{center}
\begin{tabular}{ccccc} 
\includegraphics[trim=0mm 0mm 0mm 0mm, width=.2\linewidth]
{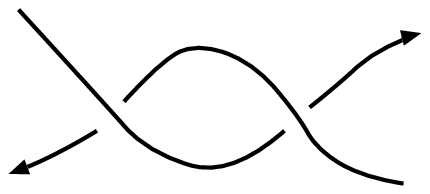},
&
\includegraphics[trim=0mm 0mm 0mm 0mm, width=.2\linewidth]
{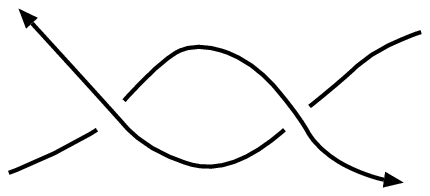},
&
\includegraphics[trim=0mm 0mm 0mm 0mm, width=.2\linewidth]
{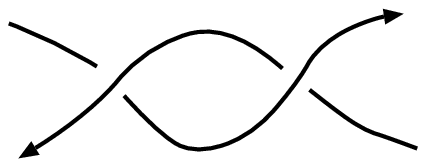}
&
or&
\includegraphics[trim=0mm 0mm 0mm 0mm, width=.2\linewidth]
{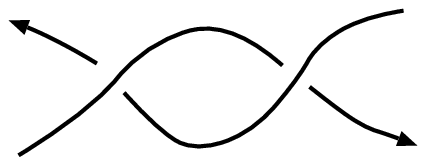}\\
\end{tabular}
\par\vs
Fig.12. Matched pairs (note the antiparallel orientations).
\end{center}
\par\vs\noindent
\par\vs
\begin{center}
\begin{tabular}{ccc} 
\includegraphics[trim=0mm 0mm 0mm 0mm, width=.28\linewidth]
{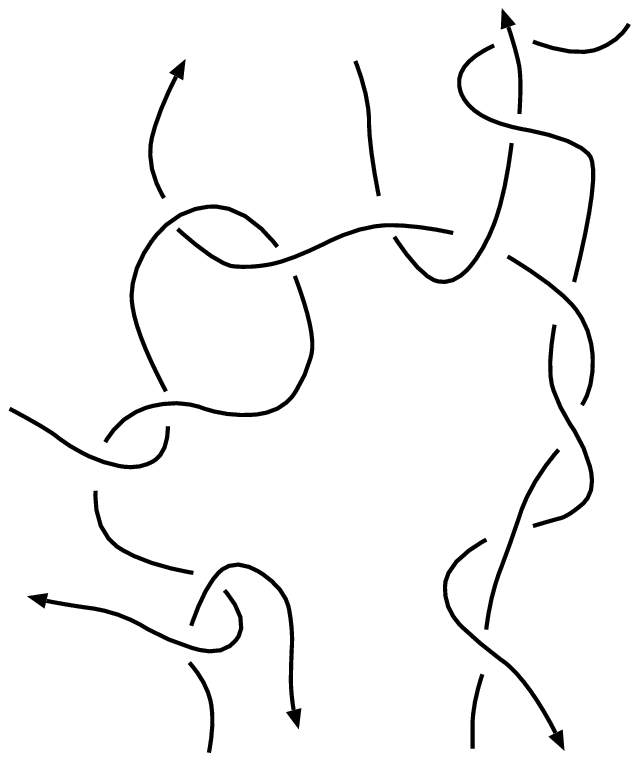}
&
\includegraphics[trim=0mm 0mm 0mm 0mm, width=.23\linewidth]
{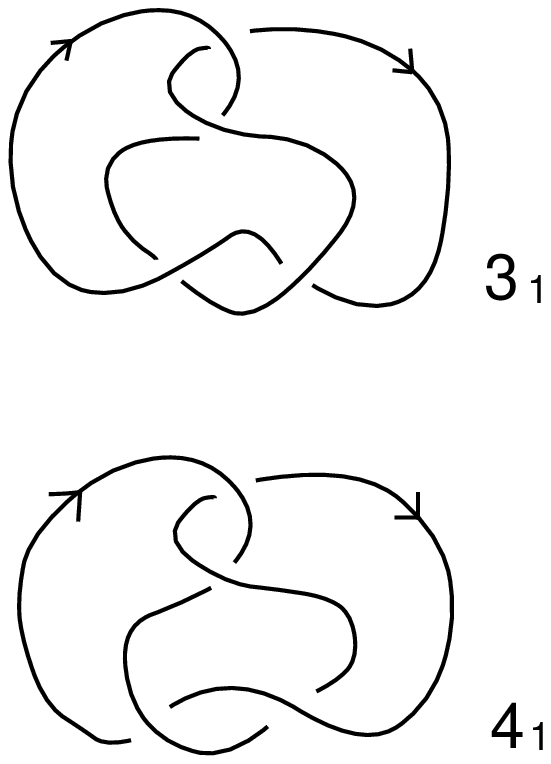}
&
\includegraphics[trim=0mm 0mm 0mm 0mm, width=.23\linewidth]
{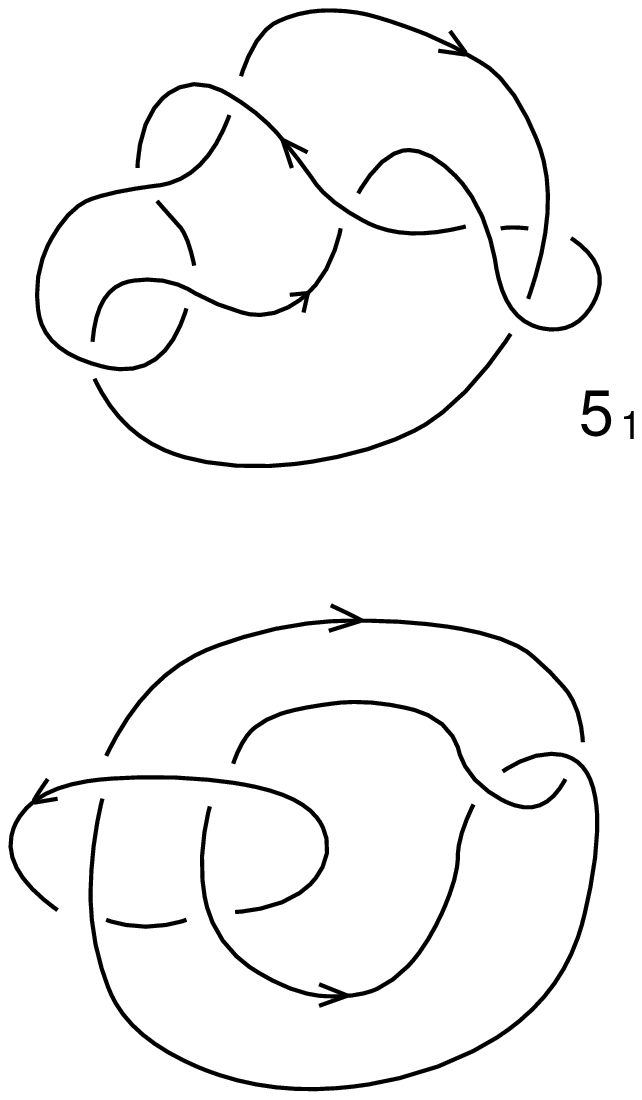}
\\
\end{tabular}
\par\vs
Fig.13. Examples of matched diagrams.
\end{center}
\par\vs\noindent
{\bf{4.8.Theorem.}} {\it{Let $L=R \cup S$ be an oriented link 
consisting of an $n$-rotor $R$ and a matched stator $S$. 
Suppose $L'$ is the corresponding $n$-rotant. Then}}
\par
(a) {\it{if $n \leq 5, L$ and $L'$ have isomorphic skein trees 
and $P_L=P_{L'}$}};
\par 
(b) {\it{if $n \leq 7$ and $R$ is cup trivial, the same conclusion holds}} ;
\par
(c) {\it{if $R$ is a parallel generalized $3$-rotor of type $2$}}
({\it{as in Theorem}} 4.4), {\it{the same}}
\par
{\it{conclusion holds}}.
\par
{\it{If $R$ and $S$ are both generalized $n$-rotors of type $k \geq 1$ and
 $S$ is matched, then $P_L(x,y,z) \equiv P_{L'}(x,y,z)$ 
modulo $n^*$ }}({\it{see Theorem}} 4.2).
\par\vs\noindent
{\bf{Proof.}} As noted in \cite{APR10}, the skein relation in 
a matched diagram looks like :
\par\vs
\begin{center}
\begin{tabular}{c} 
\includegraphics[trim=0mm 0mm 0mm 0mm, width=.5\linewidth]
{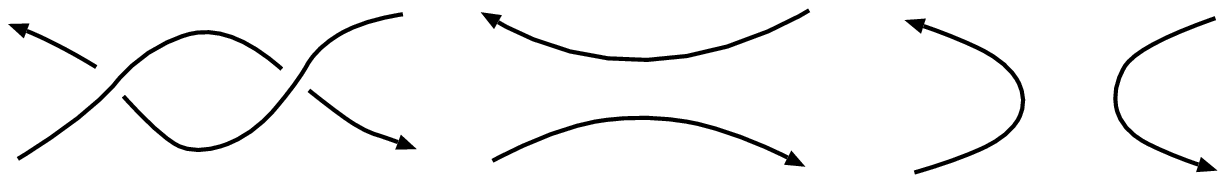}
\\
\end{tabular}
\end{center}
\par\vs\noindent
This means that one can build skein resolving trees whose leaves 
have no crossings in the stator part. Thus the arguments already 
given for the bracket resolving trees apply to the situation 
of this theorem.
\hfill$\Box$
\par\vs
In our main results (Theorem 3.1, 3.2 and 3.3) the hypotheses 
on $n$ ($\leq 5, \leq 4$ and $3$, respectively) cannot be improved. 
This is shown by examples in Jin and Rolfsen \cite{APR27}. 
Also, Example 3.6 shows the orientation assumptions are necessary.
\par\vs\noindent
{\bf{Acknowledgement}}
\par\vs
We would like to thank Steve Bleiler and Jim Hoste for useful 
conversations and Wolfgang Holtzmann, Jim Hoste and Teresa Przytycka 
for supplying computer programs for calculating link polynomials.

\par
\par\vs
%%%%%%%%%% References

\bigskip
{\small
\begin{thebibliography}{999999999}

\bibitem{APR1} R. D. Brandt, W. B. R. Lickorish and K. C. Millett, 
A polynomial invariant of unoriented knots and links, 
Invent. Math. 84 (1986) 563-573.

\bibitem{APR2} R. L. Brooks, C. A. B. Smith, A. H. Stone and W. T. Tutte, 
The dissection of rectangles into squares, Duke Math.J. 7 (1940), 312-340.

\bibitem{APR3} G. Burde and H. Zieschang, Knots (De Gruyter, Berlin, 1985).

\bibitem{APR4} J. H. Conway, An enumeration of knots and links in: 
J. Leech, ed., Computational Problems in Abstract Algebra 
(Pergamon, Oxford, 1969) 329-358.

\bibitem{APR5} R. H. Crowell and R. H. Fox, Introduction to Knot 
Theory (Ginn and Co, 1963); also Graduate Texts in Math. 
57 (Springer, Berlin).

\bibitem{APR6} P. Freyd, D. Yetter, J. Hoste, W. B. R. Lickorish,
K. Millett and A. Ocneanu, A new polynomial invariant of knots and 
links, Bull. Amer. Math. Soc. 12 (1985) 239-246.

\bibitem{APR7} C. Giller, A family of links and the Conway calculus, 
Trans. Amer. Math. Soc. 270 (1982) 75-109.

\bibitem{APR8} C. F. Ho, A new polynomial invariant for knots and links, 
preliminary report, AMS Abstracts 6 (1985) 300.

\bibitem {APR9} J. Hoste, A polynomial invariant of knots and links, 
Pacific J. Math.124 (1986) 101-108.

\bibitem {APR10} F. Jaeger, On Tutte polynomials and link polynomials, 
Proc. Amer. Math. Soc. 103 (2) (1988) 647-654.

\bibitem{APR11} V. Jones, A polynomial invariant for knots via 
Von Neumann algebras, Bull. Amer. Math. Soc. 12 (1985) 103-111.

\bibitem {APR12} L. Kauffman, State models and the Jones polynomial, 
Topology 26 (1987) 297-309.

\bibitem {APR13} L. Kauffman, On Knots, Ann. of Math. Studies 115 (1987).  

\bibitem {APR14} W. B. R. Lickorish, A relationship among link polynomials, 
Math. Proc. Cambridge Philos. Soc. 100 (1986) 109-112.

\bibitem {APR15} W. B. R. Lickorish, The panorama of polynomials for 
knots, links and skeins, Contemp. Math. 78 (1988) 399-414.

\bibitem{APR16} W. B. R. Lickorish and A. S. Lipson, Polynomials of 
2-cable-like links, Proc. Amer. Math. Soc. 100 (1987) 355-361.

\bibitem {APR17} W. B. R. Lickorish and K. Millett, A polynomial 
invariant of oriented links, Topology 26 (1987) 107-141.

\bibitem {APR18} H. R. Morton and P. Traczyk, The Jones polynomial of 
satellite knots around mutants, Contemp. Math. 78 (1988) 587-592.

\bibitem {APR19} K. Murasugi, Jones polynomials of periodic links, 
Pacific J. Math. 131 (2) (1988) 319-329.

\bibitem {APR20} J. H. Przytycki, Survey on recent invariants in classical 
knot theory, Warsaw University, Canad. J. Math. 26 (2) (1989)
\footnote{Added for e-print: It is a mistake, it should be: 
Warsaw University Preprints 6,8,9; 1986 (in English); a part of the
book: {\it Knots: a combinatorial approach to the knot theory}, (in Polish).
Script, Warsaw, August 1995, 240+ XLVIIIpp.}.

\bibitem {APR21} J. H. Przytycki, Equivalence of cables of mutants of 
knots, Canad. J. Math. 26 (2) (1989).

\bibitem {APR22} J. H. Przytycki and P. Traczyk, Invariants of links 
of Conway type, Kobe J. Math. 4 (1987) 115-139.

\bibitem {APR23} D. Rolfsen, Knots and Links (Publish or Perish Press, 
Berkeley, CA, 1976).

\bibitem {APR24} M. B. Thistlethwaite, A spanning tree expansion of 
the Jones polynomial, Topology 26 (1987) 297-309.


\bibitem {APR25} W. T. Tutte, Codichromatic graphs, J. Combin. 
Theory Ser. B 16 (1974) 168-174.

\bibitem {APR26} W. T. Tutte, Rotors in graph theory, Ann. Discrete Math. 
6 (1980) 343-347.


\bibitem {APR27} G. T. Jin and D. Rolfsen, Some remarks on rotors 
in link theory, to appear\footnote{Added for e-print: \ 
Canad. Math. Bull. 34 (1991), 480-484.}.

\end{thebibliography}}
\end{document}